\documentclass[12pt,twoside]{amsart}
\usepackage{amssymb}
\usepackage{amscd}

\title[Toric varieties]
{Toric varieties whose canonical divisors are divisible by 
their dimensions}
\author{Osamu Fujino} 
\subjclass[2000]{Primary 14M25; Secondary 14E30.}
\date{2005/1/5}
\address{Graduate School of Mathematics\\ 
 Nagoya University, Chikusa-ku Nagoya 464-8602 Japan}
\email{fujino@math.nagoya-u.ac.jp}

\newcommand{\mult}[0]{{\operatorname{mult}}}

\newtheorem{thm}{Theorem}[section]
\newtheorem{lem}[thm]{Lemma}
\newtheorem{cor}[thm]{Corollary}
\newtheorem{prop}[thm]{Proposition}
\newtheorem{claim}{Claim}

\theoremstyle{definition}

\newtheorem{rem}[thm]{Remark}
\newtheorem*{ack}{Acknowledgments}       
 
\newtheorem*{notation}{Notation}         

\begin{document}
\bibliographystyle{amsalpha+}

\begin{abstract}
We totally classify the projective toric varieties 
whose canonical divisors are divisible by their dimensions. 
In Appendix, we show that Reid's toric Mori theory implies 
Mabuchi's characterization of the projective 
space for toric varieties. 
\end{abstract}

\maketitle

\section{Introduction}
In \cite[Section 5]{hm}, 
Akio Hattori and Mikiya 
Masuda determined the structures of 
$n$-dimensional {\em{non-singular complete}} toric varieties 
whose first Chern classes are divisible by $n$ or $n+1$ as 
applications of 
their theory. Their results are as follows:  

\begin{thm}[{cf.~\cite[Corollaries 5.4, 5.8]{hm}}]\label{hatto}
Let $M$ be a {\em{complete non-singular}} toric 
variety of dimension $n$. 
\begin{itemize}
\item[(A)] If $c_1(M)$ is divisible by $n+1$, then 
$M$ is isomorphic to the projective space $\mathbb P^n$ as 
a toric variety. 
\item[(B)] If $c_1(M)$ is divisible by $n$, then 
$M$ is isomorphic to an $(n-1)$-dimensional projective space 
bundle over $\mathbb P^1$ as a toric variety. 
\end{itemize}
For 
the more precise statements, see 
\cite[Corollaries 5.6, 5.8]{hm}. 
\end{thm}

These results seem to be toric geometric analogues 
of Kobayashi-Ochiai's theorems (see \cite{ko}). In 
\cite{ko}, they characterized $n$-dimensional 
Fano manifolds whose 
first Chern classes are divisible by $n$ or $n+1$. 
Before we state the main theorem of this paper, 
let us recall the 
following theorem, which is a direct consequence 
of the main theorem of \cite{fujino1}. 

\begin{thm}\label{re00}
Let $X$ be an $n$-dimensional {\em{projective}} toric 
variety and $B=\sum_j d_jB_j$ a $\mathbb Q$-divisor on $X$, 
where $B_j$ is a torus invariant prime divisor and $0\leq d_j 
\leq 1$ for every $j$. 
Assume that $K_X+B$ is $\mathbb Q$-Cartier, 
not nef, and 
$-(K_X+B)\equiv ND$ for some Cartier divisor $D$ on $X$, 
where $N$ is a positive rational number. 
Then, \cite[Theorem 0.1]{fujino1} implies 
$N\leq n+1$. Furthermore, $N=n+1$ if and only if 
$X\simeq \mathbb P^n$, $B=0$, 
and $\mathcal O_X(D)\simeq \mathcal O_{\mathbb P^n}(1)$. 
More generally, $N>n$ implies that 
$X\simeq \mathbb P^n$, $\sum _j d_j <1$, 
and $\mathcal O_X(D)\simeq \mathcal O_{\mathbb P^n}(1)$.  
\end{thm}

Obviously, 
Theorem \ref{re00} is much stronger 
than Theorem \ref{hatto} 
(A) for {\em{projective}} toric varieties. Note that we 
do not assume that $X$ is {\em{non-singular}} 
in Theorem \ref{re00}. 
Unfortunately, we need the {\em{projectivity}} assumption 
for our 
proof since it depends on the toric Mori theory. 
In this short paper, we try to generalize 
Theorem \ref{hatto} (B) for {\em{projective}} toric 
varieties without any assumptions about singularities. 
The next theorem is the main theorem of this paper.  

\begin{thm}\label{mama}
Let $X$ be an $n$-dimensional 
{\em{projecitve}} toric variety such 
that $K_X$ is $\mathbb Q$-Cartier. 
Assume that $K_X\equiv nD$ for some Cartier divisor $D$ on 
$X$. Then, we can determine the structure of $X$. 
More precisely, if $X$ is non-singular, then $X$ has 
a $\mathbb P^{n-1}$-bundle structure over $\mathbb P^1$. 
If $X$ is singular, then $X$ is $\mathbb P(1,1,2,\cdots, 2)$ or 
the toric variety constructed in {\em{Theorem \ref{key2}}}. 
For the more precise 
statements, see {\em{Theorems {\ref{th1}}}} and 
{\em{\ref{key2}}} below. 
\end{thm}

This paper is not self-contained. It heavily relies on 
my previous paper:~\cite{fujino1}. 
As we said before, 
we need the {\em{projectivity}} assumption for 
our proof since it depends on the toric Mori theory.  
I do not know if our results are 
true or not without this assumption. 
In general, if $X$ is non-projective, then 
the Kleiman-Mori cone $\overline {NE}(X)$ may 
have little information (see \cite{klmo}). 
Finally, in Appendix, we show that Reid's toric Mori theory implies 
Mabuchi's characterization of the projective space 
for toric varieties (see Theorem \ref{mab}). 
We freely 
use the notation in \cite{fujino1}. We will work over 
an algebraically closed field $k$ throughout this note. 

We summarize the contents of this paper. 
In Section \ref{sec2}, we investigate $\mathbb Q$-factorial 
toric Fano varieties with $\rho=1$ that have long extremal 
rays. 
It is a generalization of \cite[Proposition 2.9]{fujino1}. 
Section \ref{sec3} is the main part of this paper. 
Here, we 
classify the toric varieties whose canonical divisors 
are divisible by their dimensions. Section \ref{sec4} is an 
appendix, where we treat Mabuchi's characterization of the 
projective space for toric varieties. 

\begin{ack}
I would like to thank Doctor Hiroshi Sato, who informed me of 
the paper \cite{hm}. 
\end{ack} 
\begin{notation} The symbol $\equiv$ denotes the numerical 
equivalence for $\mathbb Q$-Cartier divisors. 
\end{notation}

\section{$\mathbb Q$-factorial toric Fano varieties 
with $\rho=1$}\label{sec2}

We use the same notation as in \cite[2.8]{fujino1}. 
The following proposition is a key result in this note. 
It is a slight generalization of \cite[Proposition 2.9]{fujino1}. 
We recommend the reader to see \cite[Section 2]{fujino1} 
before reading this section. 

\begin{prop}\label{key} 
Let $X$ be an $n$-dimensional 
$\mathbb Q$-factorial toric Fano 
variety with Picard number one. 
If $X\not\simeq \mathbb P^n$ and 
$-K_{X} \cdot V(\mu_{l,m})\geq n$ for every pair $(l,m)$, 
then $X\simeq \mathbb P(1,1,2,\cdots, 2)$.  
\end{prop}
\begin{proof}
It is obvious that $n\geq 2$. 
By the assumption, we have 
$$-K_{X} \cdot V(\mu_{k,n+1})=\frac {1}{a_{n+1}}
{(\sum_{i=1}^{n+1} a_i)}
\frac{\mult {(\mu_{k,n+1})}}
{\mult {(\sigma_{k})}}\geq n
$$ 
for $1\leq k\leq n$. 
Thus 
$$
(n+1)a_{n+1}\geq \sum_{i=1}^{n+1} a_i 
\geq \frac{\mult {(\sigma_{k})}}
{\mult {(\mu_{k,n+1})}}na_{n+1}  
$$
for every $k$. 
Since 
$$
\frac{\mult {(\sigma_{k})}}
{\mult {(\mu_{k,n+1})}}\in \mathbb Z_{>0}, 
$$
we have $\mult {(\sigma_{k})}=\mult (\mu_{k,n+1})$ 
for every $k$. 
This implies that $a_{k}$ divides $a_{n+1}$ for all $k$. 
\begin{claim}\label{cl1} 
$a_1=a_2=1$, $a_3=\cdots=a_{n+1}=2$. 
\end{claim}
\begin{proof}[Proof of Claim \ref{cl1}]
If $a_1=a_{n+1}$, then $a_1=a_2=\cdots =a_{n+1}=1$ 
since we assumed $a_1\leq \cdots \leq a_{n+1}$. 
This and 
$-K_X\cdot V(\mu_{l,m})\geq n$ for every 
$(l,m)$ imply that $X\simeq \mathbb P^n$. 
See the proof of \cite
[Proposition 2.9]{fujino1}. 
Thus, we have $a_1\ne a_{n+1}$. 
It follows from this fact that $a_2\ne a_{n+1}$ since 
$v_1$ is primitive and $\sum _i a_i v_i=0$. 
In this case, 
$$-K_{X} \cdot V(\mu_{k,n+1})=\frac {1}{a_{n+1}}
{(\sum_{i=1}^{n+1} a_i)}\geq n  
$$ 
implies $a_1=a_2=1$, $a_3=\cdots=a_{n+1}=2$. 
We note that 
$$
\frac{a_i}{a_{n+1}}\leq \frac{1}{2} 
$$ 
for $i=1,2$ and $a_i\leq a_{n+1}$ for $3\leq i \leq n$. 
\end{proof}
\begin{claim}\label{cl2}  
$\mult (\sigma_1)=\mult(\sigma_2)=1$, 
that is, $\sigma_1$ and $\sigma_2$ are non-singular cones. 
\end{claim}
\begin{proof}[Proof of Claim \ref{cl2}] 
It is sufficient to prove $\mult (\sigma_1)=1$. 
We note that $\mult (\mu_{1,l})=\mult (\sigma_1)$ for 
$3\leq l\leq n+1$ and $v_2$ is primitive imply 
that all the lattice points included in 
$$
\left\{\left. \sum _{i=2}^{n+1}t_i v_i \right| 0 \leq t_i\leq 1
\right\} \subset N_{\mathbb R}
$$ 
are vertices. Thus, $\mult (\sigma_1)=1$. 
\end{proof}
Therefore, $\{v_1, v_2, \cdots, v_{n+1}\}$ spans the lattice 
$N\simeq \mathbb Z^n$. 
Thus, we obtain $X\simeq \mathbb P(1,1,2,2,\cdots, 2)$, 
a weighted projective space. 
\end{proof}

\begin{rem}\label{rem1}  
Let $X\simeq \mathbb P(1,1,2,2,\cdots, 2)$. 
Then it is not difficult to 
see that $V(v_i)$ is a torus invariant Cartier divisor 
and 
$K_X\sim -n V(v_i)$ for $3\leq i \leq n+1$. 
\end{rem}

\section{Main Theorems}\label{sec3} 

In this section, we classify the structures of the 
$\mathbb Q$-Gorenstein 
projective 
toric varieties $X$ with $-K_X\equiv nD$. 
Before we go to the classification, let us note 
the following lemma. 
The proof is easy. 

\begin{lem}[Numerical equivalence and 
$\mathbb Q$-linear equivalence]\label{lin}
Let $X$ be a projective toric variety and 
$D$ a Cartier divisor on $X$. 
Then $D\equiv 0$ if and only if $D\sim 0$. 

Let $D_1$ and $D_2$ be $\mathbb Q$-Cartier divisors on 
$X$. 
Then $D_1\equiv D_2$ if and only if $D_1
\sim_{\mathbb Q} D_2$. 
\end{lem}

First, we decide the structures of $X$ 
under the assumption 
that $X$ is $\mathbb Q$-factorial and $-K_X\equiv nD$, 
where $n=\dim X\geq 2$. 

\begin{thm}[$\mathbb Q$-factorial case]\label{th1}  
Let $X$ be a $\mathbb Q$-factorial projective toric 
variety with $\dim X=n\geq 2$. 
Let $D$ be a Cartier divisor on $X$. 
If $-K_X\equiv nD$, then the one of the following holds. 
\begin{itemize}
\item[(1)] $X\simeq \mathbb P_{\mathbb P^1}
(\mathcal O(q_1)\oplus \mathcal O (q_2)\oplus \cdots 
\oplus \mathcal O(q_n))$ such that 
$\sum _{i=1}^{n}q_i=2$. 
In this case, $\mathcal O_X(D)\simeq \mathcal O_{\mathbb P}
(1)$, 
where $\mathcal O_{\mathbb P}(1)$ is the tautological 
line bundle of 
$\mathbb P_{\mathbb P^1}
(\mathcal O(q_1)\oplus \mathcal O (q_2)\oplus \cdots 
\oplus \mathcal O(q_n))$. Note that $X$ is non-singular and 
$\rho(X)=2$. 
\item[(2)] $X\simeq \mathbb P(1,1,2,2,\cdots, 2)$, 
and $D$ is a torus invariant prime Cartier divisor on 
$X$, see {\em{Remark \ref{rem1}}}. Note 
that $X$ is singular and $\rho(X)=1$. 
\end{itemize}
\end{thm}
\begin{proof}
Since $K_X$ is not nef, there exists a $K_X$-negative 
extremal ray $R$. 
Its length is obviously $\geq n$. 
This means that $-K_X\cdot C\geq n$ for every 
integral curve $C$ such that $[C]\in R$. 
So, $\beta=0$ or $1$ in the proof of 
the theorem in \cite{fujino1} (see \cite[p.558--559]{fujino1}). 
If $\beta=1$, then it can be checked easily that 
$\alpha=0$ (see \cite[p.558]{fujino1}). 
In this case, there exists a contraction $\varphi:X\longrightarrow 
\mathbb P^1$ such that the general fibers are 
$\mathbb P^{n-1}$. 
Therefore, $F\cdot D^{n-1}=1$ for any fiber $F$ since 
$\varphi$ is flat. 
Thus, every fiber is reduced and isomorphic to 
$\mathbb P^{n-1}$. 
So, we obtain $X\simeq \mathbb P_{\mathbb P^1}
(\mathcal O(q_1)\oplus \mathcal O (q_2)\oplus \cdots 
\oplus \mathcal O(q_n))$. 
We can assume that $0<\sum q_i\leq n$ without loss of 
generality. 
Since $\mathcal O(K_X)\simeq 
\varphi^*\mathcal O_{\mathbb P^1}
(\sum q_i-2)(-n)$ and $-K_X\equiv nD$, we have 
$\sum q_i=2$. 
Therefore, $\mathcal O_X(K_X)\simeq 
\mathcal O_{\mathbb P}(-n)$. 
We finish the proof when $\beta =1$. 
When $\beta=0$, it is obvious that 
$\alpha =0$ and 
$\rho(X)=1$.  
Then this case follows from Proposition \ref{key}. 
\end{proof}

\begin{rem}\label{crepant} 
Take $X=\mathbb P_{\mathbb P^1}(\mathcal O\oplus 
\cdots\oplus\mathcal O\oplus\mathcal O(2))$, which is a 
special case of (1) in Theorem \ref{th1}. 
Then, the Picard number $\rho(X)=2$. 
So, $NE(X)$ has two rays. 
One ray $R$ corresponds to the $\mathbb P^{n-1}$-bundle 
structure $X\longrightarrow \mathbb P^1$. 
Another ray $Q$ corresponds to the contraction 
$\varphi:=\varphi_Q:X\longrightarrow \mathbb P(1,1,2,\cdots,
2)$. We note that $K_X$ is $\varphi$-numerically trivial and 
that $\varphi$ contracts a divisor $\mathbb P^1\times
\mathbb P^{n-2}\simeq \mathbb P_{\mathbb P^1}
(\mathcal O\oplus \cdots\oplus \mathcal O)\subset X$. 
Thus, $\varphi$ is a crepant resolution of $\mathbb 
P(1,1, 2, \cdots, 2)$
\end{rem}

Next, we investigate the structures of $X$ when 
$X$ is not $\mathbb Q$-factorial and $-K_X\equiv nD$. 
In the following theorem, it is obvious that 
$n\geq 3$. It is because every toric surface is 
$\mathbb Q$-factorial. 

\begin{thm}[non-$\mathbb Q$-factorial case]\label{key2} 
Let $X$ be a non-$\mathbb Q$-factorial projective toric 
variety with $\dim X=n\geq 3$. 
Assume that $X$ is 
$\mathbb Q$-Gorenstein and 
$-K_X\equiv nD$ for some Cartier divisor on $X$. 
We put $Y=\mathbb P_{\mathbb P^1}
(\mathcal O\oplus \cdots \oplus \mathcal O\oplus 
\mathcal O(1)  
\oplus \mathcal O(1))$. 
Then $X$ is the target space of the flopping contraction 
$\varphi:Y\longrightarrow X$. Note 
that $\varphi$ contracts $\mathbb P^1\times \mathbb P^{n-3}
\simeq \mathbb P_{\mathbb P^1}(\mathcal O\oplus \cdots
\oplus \mathcal O)\subset \mathbb P_{\mathbb P^1}
(\mathcal O\oplus \cdots \mathcal O\oplus 
\mathcal O(1)  
\oplus \mathcal O(1))$. 
In this case, $\rho(X)=1$ and $X$ is Gorenstein. 
\end{thm}
\begin{proof}
We take a small projective toric $\mathbb Q$-factorialization 
$f:Y\longrightarrow X$ (see \cite[Corollary 5.9]{fujino1}). 
Since $Y$ is $\mathbb Q$-factorial, 
$K_Y\equiv nf^*D$, and $\rho(Y)\geq 2$, 
we have 
$Y\simeq \mathbb P_{\mathbb P^1}
(\mathcal O(q_1)\oplus \mathcal O(q_2)\oplus \cdots 
\oplus \mathcal O(q_n))$ with 
$\sum q_i=2$. 
Since $\rho(Y)=2$, $NE(Y)$ has two rays $R$ and $Q$. 
One ray $R$ corresponds to the $\mathbb P^{n-1}$-bundle 
structure $Y\longrightarrow \mathbb P^1$. 
Another ray $Q$ corresponds to the flopping 
contraction $\varphi:=\varphi_Q:Y\longrightarrow X$. 
Note that $Q$ is spanned by one of the sections 
$C_i:=\mathbb P_{\mathbb P^1}(\mathcal O(q_i))\subset 
Y$ for $1\leq i\leq n$. It is because 
all extremal rays are spanned by torus invariant curves. 
We can assume that $q_1\leq q_2\leq \cdots \leq q_n$ 
without loss of generality. 
Since $\sum q_i=2$, we have $q_1\leq 0$. 
Note that $K_X\cdot R<0$ and $K_X\cdot C_i
=-nq_i$. 
If $q_1<0$, then $Q$ is spanned by some $C_{i_0}$ with 
$K_X\cdot C_{i_0}=-nq_{i_0}>0$. 
It is because $NE(Y)$ is spanned by $R$ and $Q$ and 
$K_X\cdot R<0$. Since $\varphi_Q$ is a flopping contraction, 
we obtain $q_1\geq 0$. 
Therefore, $(q_1, q_2, \cdots, q_n)=(0,0,\cdots,0,1,1)$ 
(see also Remark \ref{crepant}). 
It is not difficult to see that 
the target space of 
the flopping contraction $\varphi_Q:Y\longrightarrow X$ 
has the desired properties. 
\end{proof} 

\section{Appendix}\label{sec4}  
In this section, we show that Mabuchi's characterization 
of the projective space for toric varieties 
(cf.~\cite[Theorem 4.1]{mabu}) 
easily follows from \cite{reid}. 
We can skip Step 2 in the proof of \cite[Theorem 4.1]{mabu} 
by applying \cite[(2.10) Corollary]{reid}. 

\begin{thm}[{cf.~\cite[Theorem 4.1]{mabu}}]\label{mab} 
Let $V$ be an $n$-dimensional complete non-singular 
toric variety. 
Assume that the normal bundle of each torus invariant divisor is 
ample. 
Then $V\simeq \mathbb P^n$. 
\end{thm}
\begin{proof}
We note that $V$ is projective since it has ample 
line bundles. Let $\Delta$ be the 
fan corresponding to $V$. 
Take an extremal ray $R$ of $NE(V)$. 
Let $C$ be a torus invariant integral curve such that 
the numerical equivalence class of $C$ is in $R$. 
Let $\langle v_1, \cdots, v_{n-1}\rangle\in \Delta$ be the 
$(n-1)$-dimensional cone corresponding to $C$. 
Take two $n$-dimensional cones $\langle 
v_1, \cdots, v_{n-1}, v_n\rangle$ and 
$\langle 
v_1, \cdots, v_{n-1}, v_{n+1}\rangle$ from 
$\Delta$. 
Thus we have $\sum _{i=1}^{n-1}a_iv_i+v_n+v_{n-1}=0$. 
Note that $V$ is non-singular. 
We put $D_i:=V(v_i)$ for every $i$. 
Since $\mathcal O_{D_i}(D_i)$ is ample, 
we obtain that $a_i=D_1\cdots D_{i-1}\cdot D_{i}^{2}\cdot 
D_{i+1}\cdots D_{n-1}>0$ for every $i$. 
Thus, $n$-dimensional cones 
$\langle v_1, \cdots, v_{i-1}, v_{i+1}, \cdots, v_{n}, v_{n+1}
\rangle\in \Delta$ for $1\leq i\leq n-1$ (see 
\cite[(2.10) Corollary]{reid}). 
Therefore, $a_i=1$ for all $i$ since $V$ is non-singular. 
So, we obtain that $V\simeq \mathbb P^n$. 
\end{proof}

The following corollary is obvious by Theorem \ref{mab}. 
\begin{cor} 
Let $V$ be an $n$-dimensional complete non-singular 
toric variety. Then $V\simeq \mathbb P^n$ if and 
only if the tangent bundle $T_V$ is ample. 
\end{cor}
\ifx\undefined\bysame
\newcommand{\bysame|{leavemode\hbox to3em{\hrulefill}\,}
\fi

\end{document}